\documentclass[10pt,twocolumn,twoside]{IEEEtran}
\usepackage[ruled]{algorithm2e}
\usepackage{cite}
\usepackage{epstopdf}
\usepackage{graphics}
\usepackage{array}
\usepackage{color}
\usepackage{amsfonts}
\usepackage{import}
\usepackage{amsthm,amsmath}
\allowdisplaybreaks
\newtheorem{lemma}{Lemma}

\newtheorem*{remark}{Remark}

\usepackage{xpatch}
\makeatletter
\xpatchcmd{\@thm}{\thm@headpunct{.}}{\thm@headpunct{}}{}{}
\makeatother
\usepackage{multirow}
\newcolumntype{L}{>{\centering\arraybackslash}m{0.7cm}}
\graphicspath{{figures/}}
\usepackage{graphicx, pifont} 
\SetAlgorithmName{Model}{model}{list of Models}
\newcommand{\RN}[1]{%
  \textup{\uppercase\expandafter{\romannumeral#1}}%
}
\DeclareMathOperator{\Tr}{Tr}
\DeclareMathOperator{\vac}{vec}

\newcommand{\conj}[1]{\overline{#1}}
\newcommand{\norm}[1]{\left\|{#1}\right\|}

\begin{document}

\title{A Convexification Approach for Small-Signal Stability Constrained Optimal Power Flow}
%
%
%

\author{Parikshit~Pareek,~\IEEEmembership{Student~Member,~IEEE,}~and~
        Hung~D.~Nguyen$^\star$~\IEEEmembership{Member,~IEEE }\vspace{-20pt}
\thanks{$^\star$Corresponding Author}
\thanks{Authors are with School of Electrical and Electronics Engineering, Nanyang Technological University, Singapore. \textit{pare0001,hunghtd@ntu.edu.sg}}}%

\markboth{Pre-Print, Cite From IEEE TRANSACTIONS ON Control of Network Systems Vol. XX}%
{}
%



\maketitle
\begin{abstract}
In this paper, a novel convexification approach for Small-Signal Stability Constraint Optimal Power Flow (SSSC-OPF) has been presented that does not rely on eigenvalue analysis. The proposed methodology is based on the sufficient condition for the small-signal stability, developed as a Bilinear Matrix Inequality (BMI), and uses network structure-preserving Differential Algebraic Equation (DAE) modeling of the power system. The proposed formulation is based on Semi-definite Programming (SDP) and objective penalization that has been proposed for feasible solution recovery, making the method computationally efficient for large-scale systems. A vector-norm based objective penalty function has also been proposed for feasible solution recovery while working over large and dense BMIs with matrix variables. An effectiveness study carried out on WECC 9-bus, New England 39-bus, and IEEE 118-bus test systems show that the proposed method is capable of achieving a stable equilibrium point without inflicting a high stability-induced additional cost. 
\end{abstract}


\IEEEpeerreviewmaketitle

\vspace{-1em}

\section{Introduction}
\IEEEPARstart{T}{he} small-signal stability assessment is pertinent for ensuring reliable power system operation. It deals with the system's capability to maintain synchronism under the influence of small disturbances \cite{kundur1994power,song2017network}. The increasing renewable energy source integration, with uncertain and intermittent nature, has introduced issues in conventional approaches of assessing the small-signal stability \cite{quintero2014impact}. Further, the economic analysis and market-driven dispatch have become a focal point of power system operations, especially in the deregulated market environment. These factors make it difficult for an Independent Service Operator (ISO) to ensure economy and stability simultaneously. Therefore, the Small-signal Stability Constrained Optimal Power Flow (SSSC-OPF) has emerged as a tool to provide a stable and economical operating point for the power system, working over an economic objective under the technical and operational constraints \cite{chung2004generation,zarate2010opf,kodsi2007application}. The insufficiency of the damping controllers in providing small-signal stability also necessitates the SSSC-OPF formulation and solution \cite{chung2004generation}. The SSSC-OPF works in conjunction with the conventional controllers such as Power System Stabilizers (PSS) to ensure power system's small-signal stability.  

The SSSC-OPF methods proposed in the literature are focused on eigenvalue analysis and efficient computation methods for eigenvalue and its sensitivities. Some works such as \cite{chung2004generation} used numerical eigenvalue sensitivity calculations for enhancing the small-signal stability constrained power transfer capacity. The authors in \cite{zarate2010opf} used first-order Taylor approximation and considered critical eigenvalue dependencies on real power change only. The maximum singular value based stability index has been used in \cite{kodsi2007application} to provide optimal tuning to damping controllers in the electricity market. However, the linearization and modifications are performed considering the Hopf-bifurcation point \cite{zohrizadeh2018penalized} and suffer from a limited range of approximation.
The authors in \cite{condren2006expected} leverage upon the closed-form eigenvalue sensitivity formulation \cite{alvarado1999avoiding,nam2000new}, and provide an expected-security cost Optimal Power Flow (OPF) solution. Nevertheless, the requirement of matrix inversion and calculation of second-order sensitivities for Hessian makes the method time-consuming. The authors in \cite{li2013eigenvalue} used Non-Linear Semi-definite Programming (NLSDP) formulation for SSSC-OPF placing spectral abscissa constraint via a smooth nonlinear constraint. The approach does not apply to large-scale systems due to NLSDP limitations, such as dense matrix variables and involvement of the matrix inverse in obtaining small-signal stability constraints using the Lyapunov theorem.  

Recently, a gradient sampling enabled Sequential Quadratic Programming (SQP) based SSSC-OPF solution method has been proposed \cite{li2016sqp}.
Nevertheless, the nonconvex formulation requires numerical differentiation for spectral abscissa sensitivities and works over the reduced system matrix. A sequential approach for SSSC-OPF has been proposed in \cite{li2019sequential} based on the decomposition of the problem into a sequence of sub-problems. The method relies on constant updating of critical eigenvalue set due to local linearization, thus increasing problem size in each iteration to deal with the issue of critical eigenvalues becoming non-critical and vice-versa in complex plane \cite{zarate2010opf}. All these works suffer from limitations of eigenvalue analysis such as local validity, repeated computation need, and nonconvex nature of stability constraint, which make SSSC-OPF challenging to solve.

In this paper, we propose a novel convexified SSSC-OPF formulation. To the best of our knowledge, this is the first attempt to solve SSSC-OPF which does not rely on eigenvalue-based constraints for stability. The proposed method does not require prior thresholds of any stability index to impose small-signal constraints. Importantly, in this work, we use Lyapunov type stability criteria \cite{pareek2019sufficient}, represented via the Bilinear Matrix Inequality (BMI). The eigenvalue analysis is presented primarily to facilitate better understanding and presentation as the community may be more familiar with eigenvalue-based stability assessment. However, this work does not rely on such eigenvalue-based stability criteria. As a result, the proposed work is not subject to the limitations of eigenvalue and eigenvalue sensitivity approaches, i.e., local validity due to linearization and repeated calculation requirements.

The proposed sufficient Lyapunov stability criterion-based formulation works over the DAE set, thus preserving the network structure and eliminating the need to matrix inverse for calculating the reduced system matrix. A set of DAE describes both dynamic and algebraic relations arising in power systems. 

The proposed approach solves the SSSC-OPF using SDP is non-iterative fashion. Further, due to the convex nature of the formulation the proposed method is tractable for solving SSSC-OPF on large-scale systems.

The main contributions are summarized as: 
 
 \begin{enumerate}
    \item Formulation and convexification of SSSC-OPF with network structure preserving DAE, based on a sufficient condition which are applicable to most of the generator and network models.
    \item Development of an eigenvalue-analysis independent approach to handle small-signal stability constraints.
    \item Development of a novel convex relaxation and a vector-norm based SDP penalization with considerably large, dense BMI constraints in matrix variables, with feasible solution recovery.
\end{enumerate}

The structure of the proposed convexified SSSC-OPF problem can be largely divided into three stages. Firstly, in Section \ref{sec:SSSSSC-OPF}, nonconvex small-signal stability condition and OPF operational constraints are dealt and convex relaxation for both constraint types have been proposed. The stator-network couplings have been imposed and then relaxed into a convex formulation to couple the network and generator side variables. In the second stage (Section \ref{sec:penatly}), feasible solution recovery methods have been developed for proposed relaxations of stability condition Bilinear Matrix Inequality (BMI), OPF, and stator-network couplings. These recovery methods are developed as the SDP-based objective penalization functions. In the last stage (Section \ref{sec:results}), the relaxation gap and errors have been identified, calculated, and reported to be within the acceptable limits with a detailed discussion. This modular strategy is adopted to get insights into the SSSC-OPF problem while solving a convexified formulation of the same. We opt to call the proposed SSSC-OPF formulation as convexified SSSC-OPF, as it includes a strategy to convert nonconvex small-signal stability condition into a convex sufficient condition via BMI relaxation, state-of-the-art OPF relaxation, and novel stator-network equilibrium relaxation. 

\section{Convexified SSSC-OPF Problem} \label{sec:SSSSSC-OPF} 
In this paper, $\mathbb{S}^+$ denotes the set of symmetric positive semi-definite matrices, while $\mathbb{R}$ is used for indicating real numbers. $\lambda(A)$ denotes eigenvalues of real matrix $A$. The superscripts $u$ and $l$ represent upper and lower bounds of the variables, respectively. $|\cdot|$ is the absolute value operator, and $||\cdot||$ represents $Euclidean-norm$. The number of buses and generators are represented as $n_b$ and $n_g$. The real part of the largest eigenvalue of matrix $A$ is represented by $\sigma_{max}(A)$. $(\cdot)^0$ indicates OPF solution values of the variables. The complex conjugate of variables is indicated as $\conj{(\cdot)}$. Moreover, by stability, we mean small-signal stability only. $\mathbf{I}$ and $\mathbf{O}$ denote identity and null matrix of appropriate dimensions. 

In compact form, the SSSC-OPF problem is described as: 
\begin{subequations}
\begin{align} 
     \textit{min} \quad & \textit{Generation cost}\\
 \textit{s.t.} \quad  & \textit{Power\, balance \& operational \, constraints} \label{eq:SSSC-OPF_opf}\\
 & \textit{Small-signal\, stability\, constraint }\label{eq:SSSC-OPF_ss}\\
 & \textit{Stator-network\, equilibrium \,constraints } \label{eq:SSSC-OPF_eq}
\end{align}
\end{subequations}

Here, power balance and operational constraints \eqref{eq:SSSC-OPF_opf} are that of OPF problem, stability constraint \eqref{eq:SSSC-OPF_ss} is developed as a BMI, while the equilibrium constraints \eqref{eq:SSSC-OPF_eq} are enforced to establish the relation between stability and OPF variables. Each set of constraints contain non-convex relations, thus need to be convexified. The three convexifications are presented below as Convexification 1: BMI relaxation, Convexification 2: OPF relaxation, and Convexification 3: Stator-network coupling relaxation. Section \ref{sec:penatly} will duly present three corresponding feasible solution recovery schemes.

\subsection{Convex Formulation of Small-Signal Stability Constraint}

The power system is modeled using the structure-preserving nonlinear DAEs. The DAE set includes generator dynamic equations and the network algebraic relations as \cite{li2013novel}: 
\begin{subequations}\label{eq:DAE}
\begin{align}
    \mathbf{\dot x} &= f(\mathbf{x},\mathbf{y}), \\
             {0}    &= g(\mathbf{x},\mathbf{y}).
\end{align} 
\end{subequations}

Here, $\mathbf{x}\in \mathbb{R}^{n}$ is dynamic, $\mathbf{y}\in \mathbb{R}^{m}$ is algebraic variable vector respectively. The $f(\cdot)$ and $g(\cdot)$ are dynamic and algebraic equation sets, respectively. The number of dynamic and algebraic variables are represented by $n$ and $m$ respectively. The details of modelling with these equations as functions ($f$,$g$) can be obtained from \cite{nguyen2016contraction}.
 Further, we consider only equilibrium point such that for a dynamic state vector $\mathbf{x}$, there exists corresponding algebraic variable $\mathbf{y}_s(\mathbf{x})$ satisfying the algebraic constraints $g(\mathbf{x},\mathbf{y}_s(\mathbf{x})) = 0 $. In order to preserve network structure, we do not attempt to eliminate the algebraic variables as what has been done in other works such as \cite{wu1994stability}. 
The linearized DAE set can be expressed in compact form as: 
    \begin{align}\label{eq:linAE}
        E \delta \mathbf{\dot z} & = J\delta \mathbf{z}.
    \end{align}
    
Here, we define $E \in \mathbb{R}^{(n+m)\times(n+m)}$ as a diagonal matrix and $E_{ii}=1$ if $i\leq n$, else zero. Also,  we define $\mathbf{z}^T=[\mathbf{x}^T\; \mathbf{y}^T]$. $J$ is the block Jacobian matrix in (\ref{eq:linAE}) and defined as: 
\begin{align}
 J(\mathbf{z}) &=
    \begin{bmatrix}
    {\partial f}/{\partial \mathbf{x}} & {\partial f}/{\partial \mathbf{y}}\\ {\partial g}/{\partial \mathbf{x}}& {\partial g}/{\partial \mathbf{y}}  
    \end{bmatrix}=
    \begin{bmatrix}
     A & B\\ C& D
    \end{bmatrix}.\label{eq:jacobian}
\end{align}

In this work, we consider the quadratic presentation of power system DAE introduced in \cite{nguyen2016contraction}. This quadratic form leads to an affine Jacobian expressed as $J(\mathbf{z})=J_0+\sum_kJ_kz_k$ for $k=1 \dots (n+m)$, for which hereafter we use its shorthand $J$. The common practice for small-signal stability is based on necessary and sufficient condition, i.e., the system is small-signal stable if and only if all the eigenvalues of the reduced Jacobian matrix have negative real part \cite{chen1998linear}. The reduced Jacobian is obtained by eliminating algebraic variables, and has the form of $\delta\mathbf{\dot x} = J_r \, \delta\mathbf{x}$ where $J_r = A -B D^{-1}C$. This eigenvalue-based approach has a major limitation rooted in the local validity of eigenvalue and eigenvalue sensitivity. This local validity implies that repeated eigenvalue calculations are required even with the slightest of deviation in states. The reduction approach also involves matrix inversion, which introduces nonconvex, nonlinear terms, and involves a high computational cost \cite{pareek2019sufficient}. 

\subsubsection*{\textbf{Sufficient Condition for Small-Signal Stability}}
We define the Lyapunov candidate $V = \delta \mathbf{z}^T Z^T E \, \delta \mathbf{z}$. With definition of Lyapunov matrix $Z= \begin{bmatrix} P & \mathbf{O} \\ R & Q \end{bmatrix}$ and block matrix $E$, we easily obtain $V = \delta \mathbf{x}^T P \delta \mathbf{x} \succ 0$ and  with $P \succ 0$ being Lyapunov matrix. Now by differentiating $V$ along with \eqref{eq:linAE}, we get $\dot V = \delta  \mathbf{z}^T \left( J^T Z + Z^T J \right) \delta \mathbf{z}$ where we use the relation $Z^T E   = E^T Z$. A similar proof is presented in \cite{pareek2019sufficient, nguyen2016contraction}.

From Lyapunov stability criteria, the system $\delta\mathbf{\dot x} = J_r \, \delta\mathbf{x}$ is small-signal stable if and only if $\dot V \leq 0 $ \cite{chen1998linear}. This is equivalent to the following negative semidefinite ($\preceq$) relationship as:
\begin{align}\label{eq:BMIs}
   \mathbf{F}  \preceq 0.
\end{align}

Here, $\mathbf{F}=J^T Z + Z^T J$ and \eqref{eq:BMIs} is a BMI in matrix variable $Z$ and $J$. Note that this matrix $\mathbf{F}$ depends on matrix $Z$ and states $z$. In a recent work \cite{pareek2019sufficient}, the BMI condition \eqref{eq:BMIs} is investigated and is shown to be a less-conservative sufficient condition for stability.

The BMI-based stability condition offers more flexibility in terms of searching for an optimal, stable solution point in OPF. This is because BMI constraint \eqref{eq:BMIs} is expressed as functions of state variables and can be used to search for a Lyapunov matrix for a state point as $\mathbf{F}(Z,\mathbf{z})$, which is not necessarily known. Here, the control variable $\mathbf{z}$ contains the power system's variables, which are responsible for explaining the dynamic and steady-state of the system. Thus, the problem of SSSC-OPF boils down to solving the BMI constraints together with the conventional OPF constraints. However, the BMI conditions are nonlinear and nonconvex in the variables, thus bringing numerical issues when scaling. We propose a convex relaxation approach to convert the nonconvex BMI problems to convex ones for tractability purposes. 

\subsubsection*{\textbf{Convexification 1: BMI relaxation}}
The negative semi-definite (NSD) relation based BMI constraint \eqref{eq:BMIs} is nonconvex and \textit{NP-hard} to solve \cite{toker1995np}. Further, the structure of the Lyapunov matrix $Z$ makes $\mathbf{F}$ in (\ref{eq:BMIs}) dense, with matrix variables. Thus, existing methods of optimization with BMIs \cite{kheirandishfard2018convex,wang2016feasibility}, having vector variables may not be suitable  . Therefore, we propose a novel relaxation of BMI with matrix variables allowing an SDP formulation of stability constraint \eqref{eq:BMIs}.

In the following, we discuss in detail the proposed convex relaxation approach. First, to separate the nonconvexity as a quadratic term, we expand $\mathbf{F}$ as: 
\begin{align}\label{eq:expendBMI}
     J^T Z + Z^T J = (J+Z)^T(J+Z)-(Z^TZ+J^TJ).
\end{align}

This new representation can be made convex if we replace the last term with a matrix variable as used in the conventional convex relaxation approaches. In particular, we use $M$ to denote $Z^TZ + J^TJ$ and arrive at the following:
\begin{equation}
     (J+Z)^T(J+Z)-M \preceq 0. \label{eq:MBMI1}
\end{equation}

However, the nonconvex term will appear in the new matrix equality $M = Z^TZ + J^TJ$. We therefore relax this nonconvexity by imposing a new positive semi-definite (PSD) condition:
\begin{align}\label{eq:rlxM}
     M \succeq Z^TZ + J^TJ .
\end{align}

Until now, we can use the two convex sets of constraints \eqref{eq:MBMI1} and \eqref{eq:rlxM} to represent the BMI-based stability condition. However, we need to cast these new constraints in LMI form so that it can be efficiently solved using state-of-art SDP solvers for large-scale problems. The details of casting in LMI form using the Schur complement lemma \cite{boyd2004convex} will be discussed in Appendix \ref{app:1}.

Though the LMI representation can offer numerical benefits, the relaxation practice introduces gaps in the solution of the SSSC-OPF. This solution gap is discussed in the following. 

\begin{remark}\textbf{(Relaxation gap):}
The relaxation (\ref{eq:rlxM}) introduces a gap between actual nonconvex and relaxed convex feasible space, which will lead to an infeasible solution. The cause of this gap is the dissimilarity between lifting variable matrix $M$ and quadratic matrix relation $Z^TZ+J^TJ$. Further, the presence of a relaxation gap means that the solution obtained with a relaxed constraint will not be stable. Therefore, a feasible solution recovery method has been developed in Section {\ref{sec:penatly}}, for convex relaxation of dense BMIs with matrix variables. The details of the SDP penalization, as feasible solution recovery method, are presented in Section {\ref{sec:penatly}}. 
\end{remark}

 \subsection{Conventional OPF and Its SDP Relaxation}\label{subsec:relaxed OPF}
In this section, we introduce the conventional OPF and its convex relaxation. Consider a network with a set of nodes $\mathcal{N}$, generator node set $\mathcal{G} \subseteq \mathcal{N}$ and the set of branches as $\mathcal{L}$. We use $k$ as the node index. For $k$-th node, the generated power is denoted as $P_{g_k}+ j Q_{g_k}$ and the load demand is given by $P_{d_k} + j Q_{d_k}$ respectively. The node voltage is given by $V_k=V_{x_k}+jV_{y_k}$. The apparent power flow from node $k$ to $l$ is $S_{kl}=P_{kl}+jQ_{kl}$. The $p(P_{g_k})=c_{2, k}P^2_{g_k}+c_{1, k}P_{g_k}+c_{0, k}$ is quadratic cost function for real power generation with $c_{2, k},~c_{1, k}$, and $c_{0, k}$ are non-negative cost coefficients. Further, $\mathbf{Y}$ is network admittance matrix with elements ${y}_{kl}$ where $(k,l) \in \mathcal{L}$.

The OPF for the real power generation cost minimization objective can be expressed as \cite{lavaei2012zero}: 
\begin{subequations}
\begin{align}
    \textbf{minimize} & ~~~ \sum\limits_{k\in \mathcal{G}} p(P_{g_k}) \label{eq:ACOPF_obj}\\
  \textbf{subject to:}  \,   & S_{g_k}-S_{d_k}=\sum \limits_{(k,l)\in \mathcal{L}}S_{kl}\quad \forall\; k \in \mathcal{N}\label{eq:ACOPF_pbal} \\
  & V^l_{k} \leq |V_k| \leq V^u_{k} \quad \forall \; k \in \mathcal{N}\label{eq:ACOPF_vlimit}\\
    & P^l_{g_k} \leq P_{g_k} \leq P^u_{g_k} \quad \forall\; k\in \mathcal{G}\label{eq:ACOPF_pg}\\
    & Q^l_{g_k} \leq Q_{g_k} \leq Q^u_{g_k} \quad \forall\; k\in \mathcal{G}\label{eq:ACOPF_qg}\\
     & |S_{kl}| \leq S^u_{kl} \quad \forall \;(k,l) \in \mathcal{L}\label{eq:ACOPF_slimit}\\
    & S_{kl}=\conj{y}_{kl}\conj{V}_kV_k-\conj{y}_{kl}\conj{V}_kV_l \quad \forall \;(k,l) \in \mathcal{L} \label{eq:ACOPF_skl}
\end{align}
\end{subequations}

In the above OPF formulation, the superscripts $u$ or $l$ represent the upper bound and lower bound of the respective quantities. Note that we do not focus on the particular realization of the cost function as we pay more attention to the feasible space defined by the corresponding constraints. Details of development of the conventional OPF model (\ref{eq:ACOPF_obj}-\ref{eq:ACOPF_skl}) can be obtained from \cite{coffrin2016qc,lavaei2012zero,low2014convex1,low2014convex}.

\subsubsection*{\textbf{Convexification 2: OPF relaxation}}
Below we briefly discuss a popular approach from literature and present our modified ACOPF relaxation. In the ACOPF, the complexity arises in the nonconvex power flow equation \eqref{eq:ACOPF_skl}, quadratic equality having product term of voltage ($V_k \conj{V_l}$). A well-studied model replaces the voltage product term with lifting variable as \cite{lavaei2012zero,low2014convex1} :

\begin{align} \label{eq:lifting}
    W_{kl}=V_k \conj{V}_l \quad (k,l \in \mathcal{N})
\end{align}

In this work, we use the model described in \cite{lavaei2012zero} (as optimization 3) for relaxation, with lifting variable matrix $W$ and voltage variable vector $\mathbf{V}$ as: 
\begin{subequations}
\begin{align}
    \textbf{minimize} & ~~~ \sum\limits_{k\in \mathcal{G}} p(P_{g_k}) \label{eq:WACOPF_obj}\\
   \textbf{subject to:} \, & P^l_{g_k}-P_{d_k} \leq \Tr\{\mathbf{Y}_kW\} \leq P^u_{g_k}-P_{d_k} \label{eq:WACOPF_pg}\\
     & Q^l_{g_k}-Q_{d_k} \leq \Tr\{\mathbf{\overline{Y}}_kW\} \leq Q^u_{g_k}-Q_{d_k} \label{eq:WACOPF_qg}\\
    & (V^l_{k})^2 \leq \Tr \{M_kW\}\leq (V^u_{k})^2 \label{eq:WACOPF_vlimit}\\
    & \Tr\{\mathbf{Y}_{kl}W\}^2+\Tr\{\mathbf{\overline{Y}}_{kl}W\}^2 \leq (S^u_{kl})^2 \label{eq:WACOPF_slimit}\\
    & W=\mathbf{V}\mathbf{V}^T \label{eq:WACOPF_W}
\end{align}
\end{subequations}
Here, the real power generation, reactive power generation, and voltage vector $\mathbf{V}$ are:
\begin{subequations}
\begin{align}
    & P_{g_k}= \Tr\{\mathbf{Y}_kW\}+P_{d_k}, \label{eq:pgrlx}\\
    & Q_{g_k}= \Tr\{\mathbf{\overline{Y}}_kW\} +Q_{d_k}, \label{eq:qgrlx}\\
    & \mathbf{V}=[\mathbf{V}_x^T\;\; \mathbf{V}_y^T]^T. \label{eq:Vxy}
\end{align}
\end{subequations}

Here, $\mathbf{Y}_k,~\mathbf{\overline{Y}}_k,~\mathbf{Y}_{kl},~\mathbf{\overline{Y}}_{kl}$ are admittance matrices which are constructed to facilitate the relaxation of OPF problems. $M_k$ is similar to the one presented in \cite{lavaei2012zero} referring to the diagonal incident matrix and has not presented here for brevity.

Further, the definition of the voltage vector $\mathbf{V}\in \mathbb{R}^{n_b}$ in \eqref{eq:Vxy} is selected deliberately as its real and imaginary parts of voltage $V_x\,,V_y$ are control variables for stability in \eqref{eq:BMIs}. This nonconvexity of power flow constraint \eqref{eq:ACOPF_pbal} is now transferred into a single equality constraint \eqref{eq:WACOPF_W}. For relaxed OPF the equality constraint (\ref{eq:WACOPF_W}) has been replaced by: 
\begin{align}\label{eq:RCACOPF}
    W &\succeq \mathbf{V}\mathbf{V}^T.
\end{align}

The exactness of OPF has been sacrificed by replacing the nonconvex, equality constraint between $W$ and $\mathbf{V}$ (\ref{eq:WACOPF_W}) with the positive semidefinite condition (\ref{eq:RCACOPF}). The works \cite{lavaei2012zero,coffrin2016qc} discuss in detail that for various systems, the numerical solution of this OPF problem can be found with exactness as $rank(W^{o})=1$ where $W^{o}$ refers to the solution of relaxed OPF. The numerical difficulty in finding an exact rank-one solution is due to the sparse nature of the network as it creates a situation with an infinite number of solutions. Therefore, there exists a unique rank-one solution with other higher rank solutions. To avoid multiple solution issue and to improve exactness, it has been suggested to add  $10^{-5}$ per unit resistance to each ideal transformer making the system graph connected \cite{lavaei2012zero}. 

The SSSC-OPF problem in this work differs from the conventional OPF in the sense that the voltage solution needs to not only satisfy all OPF constraints but also satisfy the stability condition. In the relaxed OPF, finding the lifting variable $W$ is sufficient to find the optimal solution $\mathbf{V}$, so one only needs to impose a positive semi-definite condition on the lifting variable $W$, i.e., $W \succeq 0$. However, SSSC-OPF contains stability constraints that cannot be expressed directly in the lifting variable $W$, but the voltage vector $\mathbf{V}$. Thus, condition \eqref{eq:RCACOPF} will be used in this work instead of $W \succeq 0$. This difference is elaborated further in the following remark.

\begin{remark}
 The SDP relaxation (\ref{eq:RCACOPF}) of equality constraint (\ref{eq:WACOPF_W}) is a departure from other practices of the convex relaxation of OPF using the lifting variables. The works like \cite{lavaei2012zero,low2014convex,coffrin2016qc} propose to replace the nonconvex equality (\ref{eq:WACOPF_W}) only with a positive semi-definite constraint on $W$ ($W\succeq 0$). This is because, in the OPF problem, the voltage variable vector $\mathbf{V}$ can be replaced by $W$ entirely. Further, a \textit{rank-one} solution matrix $W^{o}$ will then uniquely decompose (eigenvalue decomposition) into voltage vector as $W^{o}=\mathbf{V}^w{\mathbf{V}^w}^T$. The voltage vector ($\mathbf{V}^w$) will then provide a feasible solution under zero-duality gap conditions \cite{lavaei2012zero}. On the contrary to the OPF, in the proposed SSSC-OPF problem formulation, $\mathbf{V}$ is a control variable for stability. Therefore, it is essential to obtain the $W^o$ such that gap between decomposed voltage vector $\mathbf{V}^w$ and control variable vector $\mathbf{V}$ is minimal if not zero. Otherwise, the OPF solution will not satisfy the stability condition and vice versa. The SDP penalization method has been proposed, in Section-\ref{sec:penatly}, to minimize the relaxation gap introduced due to PSD condition in (\ref{eq:RCACOPF}).
\end{remark}

Thus, the relaxed OPF problem is minimizing the objective \eqref{eq:WACOPF_obj} over the constraints (\ref{eq:WACOPF_pg})-(\ref{eq:WACOPF_slimit}), (\ref{eq:pgrlx}), (\ref{eq:qgrlx}) and (\ref{eq:RCACOPF}).

\subsection{Convexification of Stator-Network Equilibrium Constraints}
Unlike the conventional OPF, SSSC-OPF requires the coupling between the network and each generator. This coupling is represented by stator-network equilibrium constraints \eqref{eq:SSSC-OPF_eq} or the so-called stator-algebraic equations \cite{kundur1994power}. To the best of our knowledge, the relaxation of this equilibrium has not been proposed yet. We present the convex relaxation for such equilibrium constraints below.

Depending on the generator model, the stator-network equilibrium constraints can have a different form. In this work, we use a high order generator model-- \RN{4}-order generator model from \cite{PSAT}. We use $i$ to indicate the generator index, avoiding confusion with node index $k$. As the existing generator model includes current states, we will exclude those current quantities and replace them with voltage ones duly. In particular, by neglecting the armature resistance and leakage reactance with a steady-state condition of constant internal voltage \cite{li2016sqp}, the currents and internal voltages can be expressed as a function of $d-q$ axis stator terminal voltages ($V_{di},V_{qi}$). Thus, the following quadratic relations in voltage variables are obtained.
\begin{subequations}\label{eq:reduced_genequi}
\begin{align}
   0 &= P_{g_i}-\bigg (\frac{E_{fi}}{x_{di}}\bigg)V_{di} - \bigg (\frac{x_{di}-x_{qi}}{x_{di}x_{qi}}\bigg)V_{di} V_{qi},\\
   0 &= Q_{g_i} - \bigg (\frac{E_{fi}}{x_{di}}\bigg)V_{qi}+\bigg(\frac{1}{x_{qi}}\bigg)V_{di}^2 + \bigg (\frac{1}{x_{di}}\bigg)V_{qi}^2 .
\end{align}
\end{subequations}
Here, $x_{di}$ and $x_{qi}$ are d-axis and q-axis synchronous reactance, $E_{fi}$ is field voltage for $i-$th generator. Further, the above quadratic relations are not ready for SSSC-OPF as they are nonconvex. We further relax those as the following. 
\subsubsection*{\textbf{Convexification 3: Stator-network coupling relaxation}}
 We introduce a lifting matrix variable $W_{dq}\in \mathbb{R}^{n_g\times n_g}$. This relaxation is very similar to that of relaxed OPF and relaxed version of (\ref{eq:reduced_genequi}) can be expressed as: 
\begin{subequations}\label{eq:rlx_reduced_genequi}
\begin{align}
  P_{g_i} &= \frac{E_{fi}\,{V}_{dq_i}}{x_{di}} + \frac{(x_{di}-x_{qi})W_{dq_{i,m}}}{x_{di}x_{qi}}, \label{eq:rlx_reduced_genequi_pg}\\
   Q_{g_i} &= \dfrac{E_{fi}{V}_{dq_{m}}}{x_{di}} - \frac{W_{dq_{i,i}}}{x_{qi}}-\frac{W_{dq_{m,m}}}{x_{di}}, \label{eq:rlx_reduced_genequi_qg}\\
  W_{dq} & \succeq \mathbf{V}_{dq}\mathbf{V}^T_{dq}. \label{eq:reduced_genequi_Wdq}
\end{align}
\end{subequations}

Here, $i\in \mathcal{G}$; $m=i+n_g$; $V_{dq_i} \in \mathbf{V}_{dq}$; and $\mathbf{V}_{dq}=[\mathbf{V}^T_{d}\; \mathbf{V}^T_{q}]^T \in \mathbb{R}^{2n_g} $. The relaxed convex constraint (\ref{eq:reduced_genequi_Wdq}) has been enforced in place of the nonconvex equality  $W_{dq}=\mathbf{V}_{dq}\mathbf{V}^T_{dq}$. Therefore, a rank-one $W_{dq}$ solution ensures the equivalence between the matrix and vector variables ($W^{o}_{dq}=\mathbf{V}^{w}_{dq}{\mathbf{V}^{w}_{dq}}^T$). The reason for enforcing the semi-definite constraint as (\ref{eq:reduced_genequi_Wdq}) rather than $W_{dq}\succeq 0$ is same as explained in the context of relaxed OPF constraint (\ref{eq:RCACOPF}) as the $\mathbf{V}_{dq}$ is present in relaxed constraints (\ref{eq:rlx_reduced_genequi_pg},\,\ref{eq:rlx_reduced_genequi_qg}) and $\mathbf{V}_{dq} $ are control variables for stability in \eqref{eq:BMIs} by appearing in $\mathbf{z}$.

Note that the above stator-network equilibrium relaxation is expressed in $d-q$ form. One needs to use Park's transformation to convert those to real and imaginary nodal voltage counterparts $V_x, V_y$. The relaxation of the Park's transformation along with relaxation of the quadratic equivalence $V^2_{x_k}+V^2_{y_k}= V^2_{d_i}+ V^2_{q_i}$ for $i$-th generator is given in Appendix \ref{app:2}.

\section{Approach for Possible Feasible Solution Recovery} \label{sec:penatly}
\begin{figure}[h]
    \centering
    \vskip -1.2em
    \includegraphics[width=0.8\columnwidth]{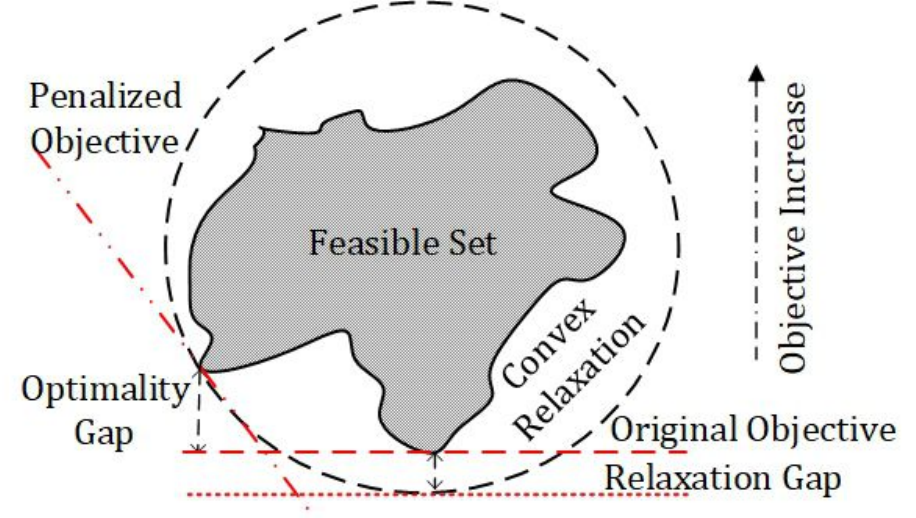}
    \vspace{-0.5em}
	\caption{Idea of objective penalization}
    \label{fig:obj_pen}
    \vspace{-1.8em}
\end{figure}

The relaxations presented in the preceding section introduce a relaxation gap, as discussed earlier. This gap will bring infeasibility and instability issues in the convexified SSSC-OPF solution. Therefore, in this section, we present a novel objective penalization approach for recovering the feasible solutions. Fig. \ref{fig:obj_pen} shows the idea behind the objective penalization in simple two-dimensional space. The dashed circle is convex relaxation of the dark shaded nonconvex feasible solution set, which introduces a relaxation gap, as shown in Fig. \ref{fig:obj_pen}. In objective penalization, the attempt is to orient the objective such that convergence will be in the feasible space. This penalization may inflict an optimality gap and hence, does not guarantee global optimality, as in Fig. \ref{fig:obj_pen}. 

\subsection{Recovery 1: Feasible BMI Solution}
As we introduce the lifting matrix variable $M$ to replace $Z^TZ+J^TJ$ in condition (\ref{eq:rlxM}), the induced relaxation gap may lead to an unstable optimal solution. One possible approach to recover the solution's feasibility is minimizing the dissimilarity between $M$ and $Z^TZ+J^TJ$. To this end, one can minimize the difference matrix $\mathbf{F}_R:=M-\{Z^TZ+J^TJ\}$ by confining its maximum (real) eigenvalue. For computational efficiency, an LMI representation of this difference matrix (details are in Appendix \ref{app:1}) is introduced as
\begin{align}\label{eq:L2}
     \mathbf{L}_2 & := \begin{bmatrix}
     M &  Z^T & J^T\\
     Z & \mathbf{I} & \mathbf{O}\\
     J & \mathbf{O} & \mathbf{I}
     \end{bmatrix}.
\end{align}

Note that, minimizing the largest eigenvalue of $\mathbf{L}_2$ will decrease that of the difference matrix $\mathbf{F}_R$ following the \textit{Interlacing theorem}. As a result, one reduces the gap between $M$ and the term $Z^TZ+J^TJ$ in \eqref{eq:rlxM}. The largest eigenvalue $\zeta$ minimization problem as SDP can be given as: 
\begin{equation}\label{eq:evp}
            \begin{aligned}
                   \min \quad & \zeta \\
                   \text{s.t.} \quad & \mathbf{L}_2-\zeta \mathbf{I} \preceq 0
            \end{aligned}
\end{equation}

Note that, this eigenvalue minimization problem (\ref{eq:evp}) introduces extra constraints where the size of $\mathbf{L}_2$ grows \textit{three times} faster than size of $J$. Also, $\mathbf{L}_2$ is dense due to $Z$ and $M$'s dense nature, which makes (\ref{eq:evp}) time consuming to solve. Another way for dissimilarity minimization can be nuclear norm $\norm{\mathbf{F_R}}_\star$ minimization, presented in  \cite{fazel2001rank}, without the condition (\ref{eq:rlxM}), an idea used extensively for matrix rank minimization. Yet, the dual formulation based LMI representation of $\norm{\mathbf{F_R}}_\star$, will be computationally complex  \cite{fazel2001rank}. Therefore, we propose below a computationally efficient alternative for imposing BMI constraints to satisfy stability conditions.

\subsubsection*{\textbf{Trace Minimization for BMI}}
Imposing the BMI condition $\mathbf{F}\preceq 0$ \eqref{eq:BMIs} will cause inefficient computation in general. A possible approach is to relax this BMI condition using its trace. The rationale behind this approach is that $\Tr\{\mathbf{F}\} \leq 0$ is a necessary condition for $\mathbf{F}\preceq 0$. 

However, minimizing the trace itself is also computationally expensive due to its nonconvex nature. We further propose to use trace's convex upper bound as the following. Note that this upper bound is valid for all BMIs of this type (\ref{eq:BMIs}).

From (\ref{eq:expendBMI}) we can obtain the trace relation as:
\begin{equation}
    \begin{aligned}
    \Tr\{\mathbf{F}\} & = \Tr \{ (Z+J)^T(Z+J)-(Z^TZ+J^TJ)\} \\
                     & = \Tr \{ (Z+J)^T(Z+J)\}-\Tr\{Z^TZ+J^TJ\}
    \end{aligned}
\end{equation}
As $Z^TZ+J^TJ \in \mathbb{S}^+$ implies $\Tr \{Z^TZ+J^TJ\} \geq 0$, thus:
\begin{equation}
    \Tr\{\mathbf{F}\} \leq \Tr \{ (Z+J)^T(Z+J)\}.
\end{equation}
Or 
\begin{equation}
\Tr\{\mathbf{F}\} \leq \norm{\vac \{Z+J\}}^2. \label{eq:trbound}
\end{equation}

The bound (\ref{eq:trbound}) is obtained using the trace-norm equality $\Tr \{ (Z+J)^T(Z+J)\} = \norm{\vac \{Z+J\}}^2$, where $\vac(X)$ indicate vectorization of matrix $X$ by arranging all columns of $X$ below each other \big(equivalent to MATLAB command $X(:)$\big). This upper bound based on $\vac(\cdot)$ is a norm-based quantity and is easily incorporated in SSSC-OPF. However, we do not impose this upper bound as constraints but construct a penalty term in the objective. The cost minimization will naturally minimize this $\norm{\vac(\cdot)}$ upper bound, thus tending to verify the BMI constraints. However, this BMI verification cannot be guaranteed due to the sufficient condition and the use of the upper estimation. An observation is that when SSSC-OPF feasible solution exists, this upper bound-based minimization leads to such a stable optimal solution. As $\norm{\vac \{\cdot\}}$ is monotonic, the objective penalty for feasible BMI solution recovery is
\begin{align}\label{eq:h1}
    h_1(Z,J)=\norm{\vac \{Z+J\}}.
\end{align}

\subsection{Recovery 2: Feasible OPF Solution}
The penalty functions for minimization of the relaxation gap in OPF will be derived from minimization of norm of difference between the variable and a base solution denoted by $(\cdot)_o$. This base solution can be a known feasible solution or the current operating point of the power system. The relaxed space encloses the actual feasible space, as shown in Fig. \ref{fig:obj_pen}. Thus, the intuitive idea is to keep the optimal solution close to the known feasible base solution. Moreover, it has been shown by authors in \cite{zohrizadeh2018penalized} that with a sufficiently large regularization parameter, similar type of objective penalization will guarantee a feasible solution with convex relaxations. This norm minimization function is
\begin{equation}\label{eq:V_Vo}
 \begin{aligned}
   ||\mathbf{V}-\mathbf{V}_o||_2^2 &= (\mathbf{V}-\mathbf{V}_o)^T(\mathbf{V}-\mathbf{V}_o)\\
           &=\mathbf{V}^T\mathbf{V}-2\mathbf{V}_o^T\mathbf{V}+\mathbf{V}_o^T\mathbf{V}_o.
\end{aligned}   
\end{equation}

The leading term of (\ref{eq:V_Vo}) is vector multiplication ($\mathbf{V}^T\mathbf{V}$) and hence cannot be minimized efficiently due to quadratic nature. To avoid this, we propose an upper bound of $\mathbf{V}^T\mathbf{V}$ by applying the trace operator on (\ref{eq:RCACOPF}) as
\begin{align}
   \Tr\{\mathbf{VV}^T\}  & \leq \Tr\{W\},\\
    \text{Or} \quad  \quad \quad \quad  \mathbf{V}^T\mathbf{V} & \leq \Tr\{W\}\label{eq:traceW}.
\end{align}

Now, we define an objective penalty function, $h_2(W,\mathbf{V})$, as upper bound of norm penalty function, $\norm{\mathbf{V}-\mathbf{V}_o}_2^2  \leq  h_2(W,\mathbf{V})$, which minimizes the gap in (\ref{eq:RCACOPF}) as

\begin{align}
       h_2(W,\mathbf{V} ) = \Tr\{W\}-2\mathbf{V}_o^T\mathbf{V}+\mathbf{V}_o^T\mathbf{V}_o. \label{eq:h2}
\end{align}

\subsection{Recovery 3: Stator-network coupling}
Similarly, the objective penalty functions for constraints (\ref{eq:reduced_genequi_Wdq}) is 
\begin{align}\label{eq:h3}
    h_3(W_{dq},\mathbf{V}_{dq}) &=\Tr\{W_{dq}\}-2\mathbf{V}_{{dq}_o}^T\mathbf{V}_{dq}+\mathbf{V}_{{dq}_o}^T\mathbf{V}_{{dq}_o}.
\end{align}

The non-convex Parks transformation is relaxed using the McCormick envelopes for bilinear terms \cite{mccormick1976computability}. We also develop penalty function $h_4$ and $h_5$ to minimize the relaxation gap imposed by Park's transformation convexification. Detailed derivation of penalties and convexification procedure is given in Appendix \ref{app:2}. The complete convexified SSSC-OPF (C-SSSC-OPF) formulation is shown as Model \ref{alg:CSSSSSC-OPF} using the weighted sum approach for handling SDP penalization with multiple penalization functions and cost minimization objective.

\begin{algorithm}[h]\label{alg:CSSSSSC-OPF}
\SetAlgoLined
\textbf{min} ~~(\ref{eq:WACOPF_obj}) + $\sum_{n=1}^{5} \gamma_n h_n$ \\
 \textbf{s.t.} ~~ (\ref{eq:WACOPF_pg})-(\ref{eq:WACOPF_slimit}), (\ref{eq:pgrlx}), (\ref{eq:qgrlx}) \\
 ~~~~~~ (\ref{eq:RCACOPF}), (\ref{eq:rlx_reduced_genequi_pg})-(\ref{eq:reduced_genequi_Wdq}), (\ref{eq:rlx_parks}), (\ref{eq:delta}), and (\ref{eq:WnWdq})\\
  \caption{Convexified SSSC-OPF (C-SSSC-OPF)}
\end{algorithm}

\section{Results and Discussion} \label{sec:results}

In this section, we present the simulation results and discussion on three different test systems. First, as in \cite{zohrizadeh2018penalized} and other works of ACOPF relaxation, we define the percentage feasibility error as the percentage difference in the trace of matrix variable and vector variable as
 \begin{align}
     \varepsilon_W & = \frac{\Tr\{W-\mathbf{V}\mathbf{V}^T\}}{\Tr\{W\}}\times 100 \% ,\label{eq:traceerror1}\\
     \varepsilon_{W_{dq}} & = \frac{\Tr\{W_{dq}-\mathbf{V}_{dq}\mathbf{V}_{dq}^T\}}{\Tr\{W_{dq}\}}\times 100 \% .\label{eq:traceerror2}
 \end{align}
 
 This error metric directly indicates the difference between the two different sets of solutions present in formulation due to the relaxation of non-linear equalities of $W = \mathbf{V}\mathbf{V}^T$ and $W_{dq}=\mathbf{V}_{dq}\mathbf{V}_{dq}^T$.  In \cite{lavaei2012zero}, the condition of zero duality is established as rank-one $W$ when only PSD condition ($W \succeq 0$) is used in the relaxation. As discussed in the remark of Section \ref{subsec:relaxed OPF}, our relaxation ($W \succeq \mathbf{V}\mathbf{V}^T$) is different from the PSD condition on $W$. Therefore, we compare the matrices by taking the ratio of their second-largest and largest eigenvalues. If eigenvalues are ordered from largest to smallest as $\lambda_1 < \lambda_2 < \dots \lambda_n$ then we define metric as
  \begin{align}
     \varepsilon_{\lambda_W} & =  \bigg \{ \frac{\lambda_2(W)}{\lambda_1(W)},\,  \frac{\lambda_2(\mathbf{V}\mathbf{V}^T)}{\lambda_1(\mathbf{V}\mathbf{V}^T))} \bigg\},\label{eq:rankerror1} \\ 
     \varepsilon_{\lambda_{W_{dq}}} & = \bigg \{ \frac{\lambda_2(W_{dq})}{\lambda_1(W_{dq})},\,  \frac{\lambda_2(\mathbf{V}_{dq}\mathbf{V}_{dq}^T)}{\lambda_1(\mathbf{V}_{dq}\mathbf{V}_{dq}^T)} \bigg\}.\label{eq:rankerror2}
  \end{align}

In case of zero-gap relaxation, $\varepsilon_{\lambda_W}$ (and $\varepsilon_{\lambda_{W_{dq}}}$) should have both entries equal to zero. In practice, we obtain small non-zero values. The target is to make the entries of $\varepsilon_{\lambda_W}$ both close to zero as possible and close to each other. Here, we want to highlight that it is important to calculate and interpret both sets of error metrics, defined in \eqref{eq:traceerror1}-\eqref{eq:traceerror2} and \eqref{eq:rankerror1}-\eqref{eq:rankerror2}, together. The percentage \textit{trace} error \eqref{eq:traceerror1} identifies the gap between two different variable sets. In contrast, the \eqref{eq:rankerror1} indicates if the penalization of objective achieves rank-one solution or not. Moreover, we do understand that the second element of $\varepsilon_{\lambda_W}$ and $\varepsilon_{\lambda_{W_{dq}}}$ will always be zero as it is the eigenvalue ratio of the rank-one matrix $\mathbf{V}\mathbf{V}^T$. However, the second element in \eqref{eq:rankerror1}-\eqref{eq:rankerror2} is reported for completeness of the results.

The mean square error (MSE) and maximum relative error (MRE) due to the convex McCormick envelopes over Park's transform (\ref{eq:rlx_parks}) is identified as $\varepsilon_{p} =\{MSE_p;\,MRE_p\}$. The $\varepsilon_{uv}= \{MSE_{uv};\,MRE_{uv}\}$ is MSE due to convexification of trigonometric equality (\ref{eq:delta}). The MRE indicates the maximum relaxation gap corresponding to a variable. At the same time, we use the MSE to analyze the effects of the objective penalties and indicate the feasibility gap. The penalties developed in section III are based on vector norms, multiplication and summation. Thus, they do not attempt to minimize the relaxation gap for individual state variables directly. This way the complicated element-wise objective is avoided. Therefore, we report both the value of MRE and the MSE related to trigonometric equality and Park's transform relaxation. The change in cost of generation in convexified SSSC-OPF (C-SSSC-OPF), as compared to that of the relaxed ACOPF solution, is indicated as $\Delta \text{Cost}\%$. The real-part of an eigenvalue is indicated as $\sigma$ with $\sigma_{max}$ being the largest real-part of eigenvalues. The PSS devices are not considered for testing and validating the proposed method's capability in enforcing small-signal stability. However, PSS can be engaged upon obtaining optimal stable solutions to stabilize the system in the presence of unexpected disturbances and variations. We model the automatic voltage regulator (AVR) via its dynamic equations in the contraction of matrix $J$. The IEEE Model I (AVR Type II of PSAT \cite{PSAT}) is followed in this work. During the small-signal stability study, the AVR control operation is considered constant \cite{li2016sqp,li2019sequential}.

\subsection{WECC 9-bus, 3-Machine System}

  \begin{table*}[t]
  \centering
  \caption{WSCC 9-bus, 3-Machine System}
          \bgroup
\def\arraystretch{1.2}\textbf{}
  \vskip -2em
    \begin{tabular}{c|c|ccccccc}
     & $\sigma_{max}$ & Cost $\$/hr$ & $\varepsilon_W$ & $\varepsilon_{W_{dq}}$ &  $\varepsilon_{\lambda_W}$ & $\varepsilon_{\lambda_{W_{dq}}}$ & $\varepsilon_{uv}$ & $\varepsilon_p$\\
                        \hline
   C-SSSC-OPF & -0.2794 & 5508.23 & 4E-08\% & 8E-09\% & \{1E-11, 0.0\} & \{5E-12, 0.0\} &\{3E-20; 5.8E-10\} & \{0.011; 0.10\}\\
    \hline
Relaxed ACOPF  & 8.915 & 5324.32  & 6E-12\% & - & \{5E-14, 0.0\} & \{-,-\} & \{-,-\} & \{-,-\} \\ 
    \hline
    \end{tabular}%
    \egroup
    \vskip -1.5em
  \label{tab:9Bus}%
\end{table*}%

\begin{table}[b]
  \centering
  \vspace{-1.3em}
  \caption{Generator set point in $pu$ for 9-bus System}
  \vspace{-1em}
    \bgroup
\def\arraystretch{1.25}
\begin{tabular}{c|ccc}
     & $P_{g_1}$ &  $P_{g_2}$ & $P_{g_{3}}$  \\
          \hline
    Relaxed ACOPF & 0.8723  & 1.3940 & 0.9273 \\
    C-SSSC-OPF & 1.3760  & 1.2330 & 0.6033 \\
    \hline
    \end{tabular}
       \egroup
  \label{tab:9buspg}
\end{table}%

\begin{figure}[t]
    \centering
    \includegraphics[width=0.8\columnwidth]{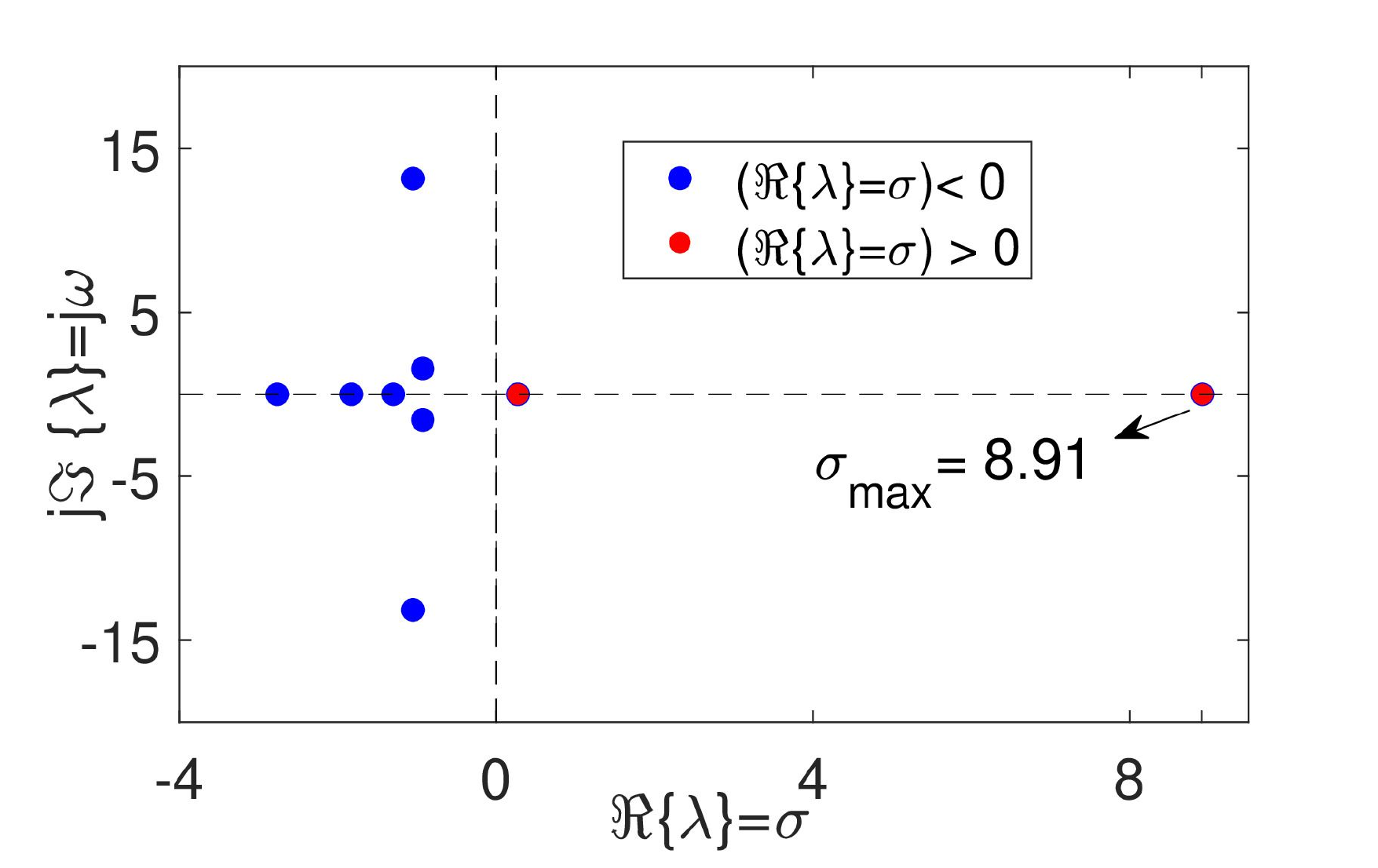}
    \vspace{-1em}
	\caption{Eigenvalue subset of relaxed OPF solution for 9-bus system}
    \label{fig:eig_opf9}
    \vspace{-1.1em}
\end{figure}

The dynamic data of this system has been taken from \cite{sauer1998power}, while the generation cost function coefficients are obtained from MATPOWER \cite{zimmerman2010matpower}. The comparative cost analysis of the proposed SSSC-OPF solution has been performed with a relaxed ACOPF solution obtained by solving relaxed ACOPF \eqref{eq:WACOPF_obj}-\eqref{eq:WACOPF_slimit},\eqref{eq:RCACOPF}. 
Fig. \ref{fig:eig_opf9} shows a set of eigenvalues with largest real-part $\Re\{\lambda\}=\sigma$, in complex plane. As indicated in the plot, the $\sigma_{max}$ is $8.91$, and the system is unstable, having two eigenvalues on the right side of $j\omega-axis$. The multiple eigenvalues on the right half of the plane present challenges to the methods of stability recovery. It is difficult to move eigenvalues while satisfying the ACOPF constraints and maintaining optimality. Therefore, the SSSC-OPF solution is essential as it will provide an economical and stable operating point.

\begin{figure}[t]
    \centering
    \includegraphics[width=\columnwidth]{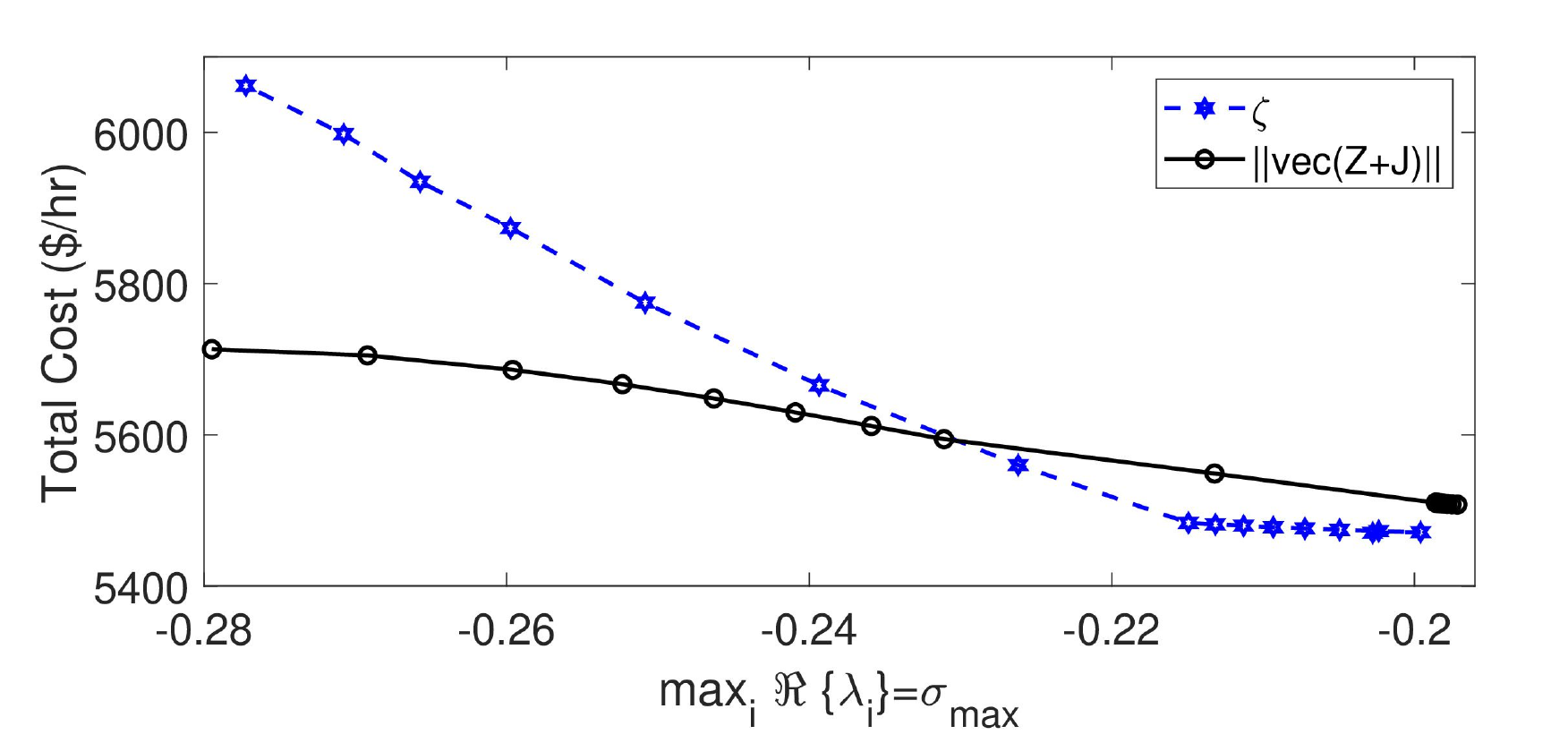}
    \vspace{-1.5em}
	\caption{Total generation cost variation with $\sigma_{max}$ for 9-bus system}
    \label{fig:cost9}
          \vspace{-1.5em}
\end{figure}

For 9-bus system, we have solved the convexified SSSC-OPF using both approaches of possible feasible BMI solution recovery, $\zeta$ with (\ref{eq:evp}) and $h_1=||\vac\{Z+J\}||$. The variation of the total cost of real power generation  with the maximum real-part of critical eigenvalue is shown in Fig. \ref{fig:cost9}. This figure shows that the proposed SDP penalization methods have been able to achieve an optimal stable solution in situations where OPF solution is not stable itself (Fig. \ref{fig:eig_opf9}). The BMI trace minimization based penalization function performs better as we move to the lower values of $\sigma_{max}$. The cost of stability is less with  $||\vac\{Z+J\}||$ than that of using the largest eigenvalue $\zeta$ penalization \eqref{eq:evp}. 

It is relevant to note that the spectral abscissa ($\sigma_{max}$) is a nonconvex function, and its relationship is difficult to identify analytically with  $||\vac(Z+J)||$ due to the involvement of Lyapunov matrix $Z$. Further, the SDP penalization function $||\vac(Z+J)||$ as well as $\zeta$ are not designed to minimize the $\sigma_{max}$ but to converge on a stable solution with $\sigma_{max} < 0$. Therefore, the results show that the proposed SSSC-OPF solution method achieves its desired objective, and the cost variation is just indicative that total generation cost increases with a decrease in $\sigma_{max}$. 

The solution of C-SSSC-OPF must be interpreted while considering all the errors present due to the relaxation gap simultaneously. Table \ref{tab:9Bus} represents C-SSSC-OPF solutions for the 9-bus system with relaxed ACOPF results. The results show that error $\varepsilon_p$ is the most significant factor, while others are negligible. Furthermore, the results confirm that, the rank-approximation errors remain close to zero for both lifting variables $W, W_{dq}$. For decreasing $\varepsilon_p$, tighter convex envelopes can be applied. Here, we highlight that critical eigenvalues obtained using vector variables and decomposed variables (from $W, W_{dq}, \mathbf{U}_u, \mathbf{U}_v$) have a gap in the order of $10^{-15}$ in all cases, indicating that error values are within acceptable limits.

The real power set-points for generators are given in Table \ref{tab:9buspg}, for an unstable relaxed ACOPF and stable C-SSSC-OPF solution. An important observation is that a higher inertia generator, $1^{st}$ generator, shares a large section of power in the C-SSSC-OPF solution compared to low inertia generators ($3^{rd}$). This observation is in line with the understanding that higher inertia generators should share more power for stability. It is worth noting here that real power-sharing will also be influenced by the cost curve, especially under high loading.

\begin{figure}[b]
    \centering
    \vspace{-1.5em}
    \includegraphics[width=0.8\columnwidth]{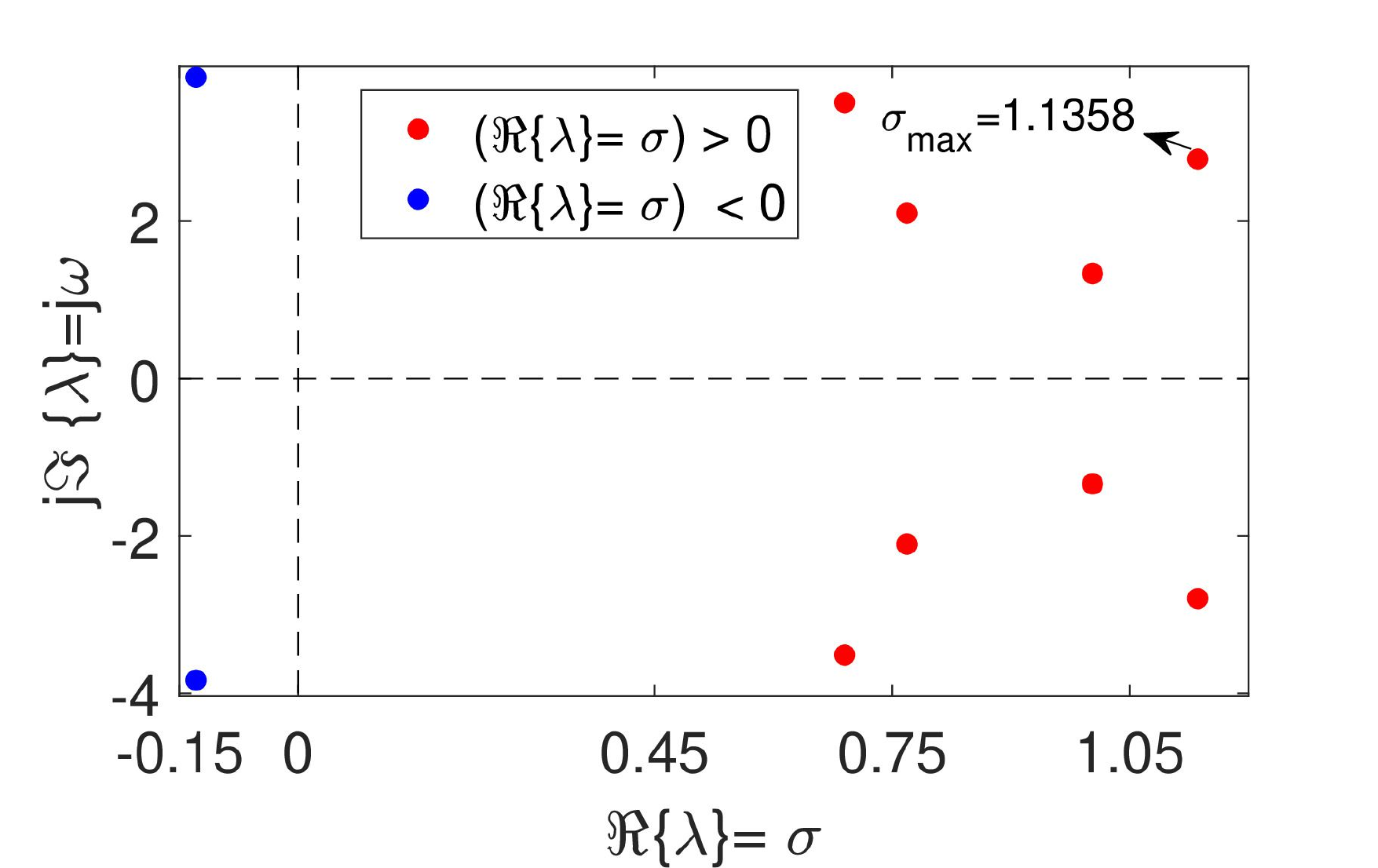}
    \vspace{-0.6em}
	\caption{Eigenvalue with largest $\sigma$ at relaxed OPF solution for 39-bus System}
    \label{fig:eig_opf39}
\end{figure}

\begin{table*}[h]
  \centering
  \caption{New England 39-bus, Ten-Machine System}
          \bgroup
\def\arraystretch{1.2}
  \vskip -1.2em
    \begin{tabular}{c|c|ccccccc}
     & $\sigma_{max}$ & Cost $\$/hr$ & $\varepsilon_W$ & $\varepsilon_{W_{dq}}$ &  $\varepsilon_{\lambda_W}$ & $\varepsilon_{\lambda_{W_{dq}}}$ & $\varepsilon_{uv}$ & $\varepsilon_{p}$\\
                        \hline
   C-SSSC-OPF & -0.1395 & 42823.7 & 3E-07\% & 0.62\% & \{5E-09, 0.0\} & \{0.0059,0.00\} & \{0.016; 0.018\} & \{0.006; 0.0935\}\\
    \hline
Relaxed ACOPF  & 1.1358 & 40951.7  & 2E-09\% & - & \{3E-12., 0.0\} & \{-,-\} & - & - \\  
    \hline
    \end{tabular}%
    \egroup
      \vspace{-1.5em}
  \label{tab:39 Bus}%
\end{table*}%

\subsection{New England 39-bus, 10-Machine System}
The system dynamic data are obtained from \cite{ramos2015benchmark}, cost coefficients from \cite{zimmerman2010matpower}, and generator limits from \cite{zarate2010opf}. For the C-SSSC-OPF solution of this system, we only use $||\vac\{Z+J\}||$ penalization function due to its computational advantages. Fig. \ref{fig:eig_opf39} shows a set of ten critical eigenvalues at relaxed ACOPF solution point, with many being on the right-hand side of $j\omega$-axis representing an instability solution. It is difficult to recover a stable solution from this situation, as described in the 9-bus test system results. 

\begin{table}[b]
  \centering
  \vskip -1.5em
  \caption{Generator set point in $pu.$ for 39-bus System}
  \vspace{-1em}
      \bgroup
\def\arraystretch{1.25}
\begin{tabular}{c|ccc}
     & $P_{g_1}$&  $P_{g_2}$ & $P_{g_{10}}$  \\
          \hline
    Relaxed OPF & \textbf{4.025}  & 7.348 & 8.558 \\
    C-SSSC-OPF & \textbf{4.025}  & \textbf{7.475} & 10.920 \\
    \hline
    \end{tabular}
    \egroup
  \label{tab:39buspg}
\end{table}%

\begin{figure}[t]
    \centering
    \includegraphics[width=\columnwidth]{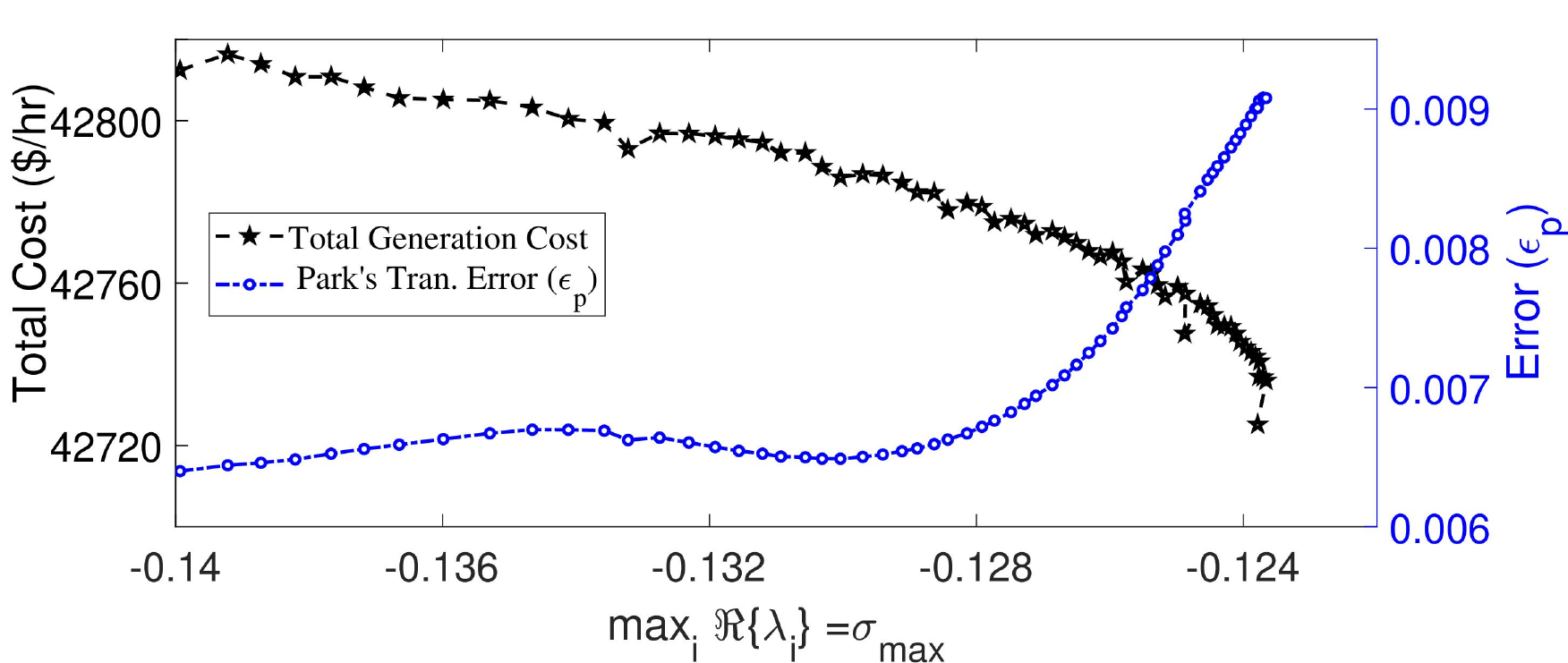}
    \vspace{-2em}
	\caption{Total generation cost and $\epsilon_p$ variations with $\sigma_{max}$ for 39-bus system}
    \label{fig:cost39}
          \vspace{-1.2em}
\end{figure}
\begin{figure}[t]
    \centering
    \includegraphics[width=\columnwidth]{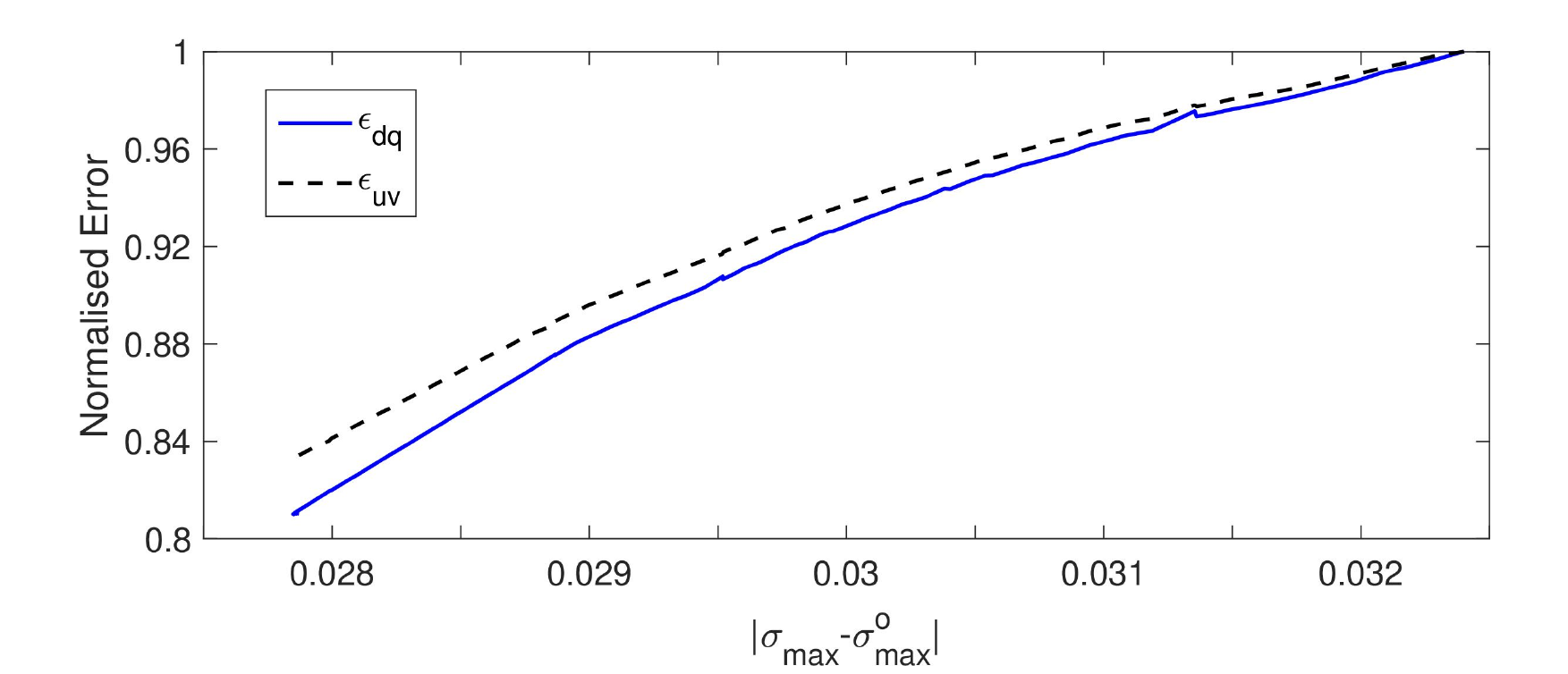}
    \vspace{-1.9em}
	\caption{Variation of absolute value of difference between $\sigma_{max}$ ans $\sigma^{o}_{max}$ with normalised error values for 39-bus system. The $\sigma^{o}_{max}$ is obtained by using the decomposed OPF variables vectors $\mathbf{V}^w$, $\mathbf{V}^w_{dq}$, $\sqrt{\mathbf{U}}$ and $\sqrt{\mathbf{V}}$.}
    \label{fig:sigma_error_39}
          \vspace{-1.5em}
\end{figure}

Fig. \ref{fig:cost39} depicts variation in optimal generation cost and $\varepsilon_p$ with $\sigma_{max}$ obtained as proposed SSSC-OPF solution. The general trend is that the cost of generation decreases with increase in $\sigma_{max}$, system getting less stable. Further, both these variations also indicate the nonconvex nature of the relationship of $\sigma_{max}$ with C-SSSC-OPF objective and error values. Table \ref{tab:39 Bus} shows the results obtained by solving the proposed convexified SSSC-OPF. First, an important point to note is that proposed method has been able to find a stable optimal solution with 4.57\% additional cost, for a system which is not stable at base ACOPF solution. The feasibility gap metric reveals that the proposed method has obtained solution with $\varepsilon_{\lambda_{W_{dq}}} < 1\%$ (relaxation gap in stator-network equilibrium) while ACOPF relaxation gap $\varepsilon_{\lambda_{W}}$ is close to zero. The results also call for the need for evaluation of effects of Park's transform relaxation on stability. Thus, Fig. \ref{fig:sigma_error_39} shows the gap between two $\sigma_{max}$ values obtained from the different sets of control variables of C-SSSC-OPF. The normalization is done based upon the maximum value of respective errors. The trend supplements our understanding that error influences the stability, and the gap increases as the error values move towards their respective maximum values. An insight from this result is that $\sigma_{max}$ value should be away from $j\omega$ axis to keep the system stable with an error. Further, the trade-off between error and $\sigma_{max}$ should be taken with consideration of the corresponding gap imposed by the errors on that particular $\sigma_{max}$. In simulations, we obtain an optimal stable solution from both $\sigma_{max}$ and $\sigma^0_{max}$.
 
 The values of $P_{g_i}$ for three generators has been presented in Table \ref{tab:39buspg}. The generator 10, which is an equivalent representation of the New York network with large inertia, shares the largest power set-point change occurring from OPF to C-SSSC-OPF. The generators 1 and 2 having large inertia are set on their maximum limits. To evaluate the effect of the cost curve on $P_{g_i}$'s we attempted a solution with increasing the $P^u_{g_4}$ value from $4.025\, pu$ to $5.00\, pu$. The proposed SSSC-OPF converges at $P_{G_4}=5.00\, pu$ with an equivalent decrease in $P_{G_{10}}$. This indicates that in the proposed SSSC-OPF formulation, cost minimization objective is influencing the generator set points along with the inertia values at higher loading. 

\begin{table*}[h]
  \centering
  \caption{IEEE 118-Bus System, 54 Machine}
            \bgroup
\def\arraystretch{1.2}
  \vskip -1.2em
    \begin{tabular}{c|c|ccccccc}
    & $\sigma_{max}$ & Cost $\$/hr$ & $\varepsilon_W$ & $\varepsilon_{W_{dq}}$ &  $\varepsilon_{\lambda_W}$ & $\varepsilon_{\lambda_{W_{dq}}}$ & $\varepsilon_{uv}$ & $\varepsilon_{p}$ \\
                        \hline
    C-SSSC-OPF &  -0.0507 & 131,041  & 0.08\% & 0.63\% & \{ 8E-4,\,0.0\} & \{0.001,\,0.000\} & \{1E-08; 7.3E-05\} & \{0.0031; 0.047\}  \\
    \hline
   Relaxed ACOPF & 0.0463 & 129,632  & 0.13\% & - & \{1E-3,\,0.0\} & \{-,\,-\} & - & -\\
    \hline
    \end{tabular}%
    \egroup
  \label{tab:118 Bus}%
  \vspace{-2em}
\end{table*}%

\subsection{IEEE 118-bus, 54-Machine System}
For 118-Bus system, we obtain the steady-state load, cost and network data from \cite{zimmerman2010matpower}. The dynamic system parameters are taken from \cite{data118}. We use the norm penalization method $h_1(Z,J)$ from \eqref{eq:h1}, instead of using largest eigenvalue penalization \eqref{eq:evp}, to make finding stable solution more tractable. 

Now from results shown in Table \ref{tab:118 Bus}, we can see that the proposed C-SSSC-OPF has been able to achieve the stable solution with additional cost of 1.09\% compared to relaxed ACOPF result. Further, the feasibility gap is 0.08\% and 0.63\% in ACOPF and stator-network equilibrium relaxation respectively. The eigenvalue metrics $\varepsilon_{\lambda_W}$ and $\varepsilon_{\lambda_{W_{dq}}}$ also suggest that the gap between lifting variable and original variable is of order of $10^{-3}$. Importantly, the eigenvalue ratio for $W$ and $W_{dq}$ is close to zero. This implies that the proposed method has been able to achieve approximately rank-one matrix variables. Furthermore, we obtain the MSE error values $0.0031$ for Park's transform and $10^{-8}$ for trigonometric equality. Like the other systems, we obtain the MRE values approximately one order of magnitude higher than MSE values in Table \ref{tab:118 Bus}. Moreover, the error in Park's transformation is more dominant than the error in trigonometric inequality for both MSE and MRE metrics. This highlights the need to develop tighter envelope relaxation of Park's transformation, which will be explored further in future works. Another important feature, highlighted by these results, is that proposed SSSC-OPF method has been able to obtain stable solution even when base ACOPF solution is not stable. Below we briefly outline the method of obtaining the weights $\gamma_n$ and numerical values of the same.

\subsubsection*{Numerical Values of $\gamma_n$}
The calculation of the weights for a multi-objective optimization problem, such as Model 1, is a non-trivial problem. However, we can estimate the weights using information about the different objectives and then improve them with a few successive iterations. The first consideration is that different objective functions \eqref{eq:WACOPF_obj} and $h_1 \dots h_5$ will be effectively minimized if they are of the same order approximately upon multiplication of weights $\gamma_n$. By ACOPF results, we know the order of generation cost for different systems. The weights $\gamma_n$ should be such that $\gamma_n h_n$ are of same order approximately. It is easy to identify that the objective penalties $h_2$ to $h_5$ have value in single digit and mostly less than one. Therefore, we start with $\gamma_n$ values, which brings penalties $h_2 \dots h_5$ to the order of generation cost. Now, a lesser value of the weight $\gamma_1$ is needed as $h_1$ has the norm of vector form of Jacobian and a variable matrix $Z$, leading to higher numerical value itself. Further, in future works, we will attempt to understand the behavior of generation cost, stability, and feasibility gap in relationship with these weights. The weights used for different systems are given in Table \ref{tab:gamma}.

\begin{table}[h]
  \centering
  \caption{Numerical values of penalty weights $\gamma_n$}
            \bgroup
\def\arraystretch{1.2}
    \begin{tabular}{c|ccccc}
   System & $\gamma_1$ & $\gamma_2$ & $\gamma_3$ & $\gamma_4$ & $\gamma_5$ \\
    \hline
    9-Bus &  1 & 500  & 1E03 & 1E03 & 1E03 \\
    \hline
   39-Bus &  10 & 2E04  & 1E04 & 1E04 & 1E04 \\
    \hline
   118-Bus &  1 & 2E04  & 3E04 & 6E04 & 6E04 \\
    \hline
    \end{tabular}%
    \egroup
  \label{tab:gamma}%
\end{table}%

In the next subsection, we present a comparative analysis with different existing works. The analysis focuses on proposed method's ability of obtaining a stable solution with less additional cost over ACOPF cost. We also discuss the computation time results, with ideas to improve the performance further. 

\subsection{Comparative Results}
Firstly, in works such as \cite{li2019sequential,li2016sqp,li2013eigenvalue}, majority of the simulations are performed on systems with base ACOPF solutions that are already small-signal stable. Further, not all works use the same definition of stability as that in the proposed work. The works \cite{li2016sqp,li2013eigenvalue} and the proposed work evaluate stability at a solution point using largest real-part of eigenvalue of Jacobian, $\sigma_{max}\leq 0$. In contrast, \cite{li2019sequential} uses minimum damping ratio (MDR) based definition and only tracks a small set of eigenvalue pairs to obtain highest damping ratio solution. The work is also focused on enhancement of the stability, not on obtaining the stable solution. Furthermore, all existing works used very different computation machines, having large variations in computation power. Considering all these factors, we attempt to compare the proposed method with the most suitable information available in these works for different test systems below. 

An interesting point to note from results in Table \ref{tab:Comp_9 Busa}-\ref{tab:Comp_118 Busa} is the increase in cost with decrease in the value of $\sigma_{max}$. The comparison is particularly important as different works use different cost functions. For that we present values of change in cost, $\Delta \text{cost}\%$, and distance on real-axis between largest real-part of eigenvalues obtained with base ACOPF and SSSC-OPF solution as $\Delta \sigma_{max}$. Although we understand that the variation in cost with variation in $\sigma_{max}$ is not linear, yet large  $\Delta \sigma_{max}$ with lower $\Delta \text{cost}\%$ does indicate the ability of a method to obtain stable solution with least cost. As indicated in table \ref{tab:Comp_9 Busa}, \ref{tab:Comp_39 Busa}, and \ref{tab:Comp_118 Busa}, the proposed method has been able to move $\sigma_{max}$ largest for every percentage increase in cost. Therefore, it shows the proposed method's ability to obtain stable solution with least cost of stability. 

\begin{table}[t]
  \centering
  \caption{Comparison Results for WSCC 9-bus, 3-Machine System}
   \bgroup
     \def\arraystretch{1.2}
  \vskip -1.0em
    \begin{tabular}{c|c|ccccc}
     & $\sigma_{max}$& $\Delta \sigma_{max}$ & $\Delta$Cost$\%$&  ACOPF &  Time(s) \\
                        \hline
  Proposed & -0.279 & 9.189 & 3.46 &  Unstable  & $<$ 1   \\
    \hline
NLSDP \cite{li2013eigenvalue}& -0.525  & 0.325 & 3.60 &  Stable & ~60  \\
\hline
SQP-GS \cite{li2016sqp}  & -0.500 & 0.460 & 10.91  & Stable  & - \\  
    \hline
    \end{tabular}%
    \egroup
    \vskip -1.5em
  \label{tab:Comp_9 Busa}%
\end{table}%

\begin{table}[t]
  \centering
  \caption{Comparison Results for 39-bus, 10-Machine System}
    \bgroup
     \def\arraystretch{1.2}
  \vskip -1.0em
    \begin{tabular}{c|c|ccccc}
     & $\sigma_{max}$& $\Delta \sigma_{max}$ & $\Delta$Cost $\%$ & ACOPF  &  Time(s) \\
                        \hline
  Proposed & -0.139 & 1.414 & 4.57 & Unstable  & 31.70   \\
    \hline
SQP-GS \cite{li2016sqp} & -0.20 & 0.09 & 3.17 & Stable  & 40.60  \\
\hline
    \end{tabular}%
    \egroup
    \vskip -1.5em
  \label{tab:Comp_39 Busa}%
\end{table}%

\begin{table}[t]
  \centering
  \caption{Comparison Results for IEEE 118-bus, 54-Machine System}
     \bgroup
     \def\arraystretch{1.2}
  \vskip -1.0em
    \begin{tabular}{c|c|ccccc}
     & $\sigma_{max}$& $\Delta \sigma_{max}$ & $\Delta$Cost$\%$ & ACOPF &  Time(s) \\
                        \hline
  Proposed & -0.0507 & 0.097 & 1.09 & Unstable  & 140.54  \\
    \hline
SQP-GS \cite{li2016sqp} & -0.10 & 0.45 & 5.73 & Unstable  & 114.93  \\
\hline
    \end{tabular}%
    \egroup
    \vskip -1.5em
  \label{tab:Comp_118 Busa}%
\end{table}%

The MDR-based definition of small-signal stability is used by authors in \cite{li2019sequential} for stability enhancement. The results presented in \cite{li2019sequential} for cost analysis are very different than obtained using the proposed method. In \cite{li2019sequential}, all the case studies are presented for stability enhancement with systems being stable at base ACOPF solution without any stability constraint imposition. Thus, the results are for cost of stability enhancement via increasing MDR, not of obtaining stability for an unstable base point. In comparison, as indicated before, the proposed work demonstrates the capability of the method to obtain stable solution when system is not stable as base ACOPF point. For 39-Bus system, as observed from Fig. 6 and Fig.7 of \cite{li2019sequential}, the $\sigma_{max}\approx-0.32$, with starting from $\sigma_{max}\approx -0.1$, and the cost change is $4.46\%$. The total time taken in solving the problem is 86.39 sec.

As for elapsed time comparisons, we compare the computation time required by the proposed method, implemented with unoptimized codes, with that reported in other works using other methods. As indicated in Table \ref{tab:Comp_9 Busa}-\ref{tab:Comp_118 Busa} the proposed method performs as good as the existing methods of SSSC-OPF. Further, as mentioned before, the difference in computation power used to solve problems is significant. In Table \ref{tab:machinea}, we list the different computer configurations provided in different works. It can be seen that the state-of-art work of SSSC-OPF presented in \cite{li2016sqp} uses much more powerful processing unit. Nevertheless, the proposed method performs relatively the same in terms of the computation time. 

Now, below we suggest a simple solution for improving the computation time performance. We refer the computation time as the sum of the elapsed time taken by YALMIP April-2019 version \cite{yalmip} and CVX version 3 \cite{cvx} for composing and the elapsed time taken by the solver to find the optimal solutions. We use MOSEK 9.2 \cite{mosek} to solve the convex optimization problems. Thus, the computation time required by proposed method can easily be reduced by using more efficient commercial solvers such as Mosek. For example, out of total time of approximately 140 seconds needed to solve the proposed C-SSSC-OPF for 118-Bus system, the MOSEK 9.2 solver only takes about 3 seconds. The rest of the elapsed time goes in the composing the problem within CVX \cite{cvx}. Therefore, there is a scope of improvement in terms of the computation time with more efficient use of solvers, which will be explored in future works.    

\begin{table}[t]
    \caption{Computing Resource used in \cite{li2016sqp} and Proposed Work}
    \centering
    \vskip -1.1em
      \bgroup
\def\arraystretch{1.2}
    \begin{tabular}{c|c|c|c}
         & Machine & Processor & RAM \\
         \hline
     Proposed    &  Intel Xenon & 3.7 GHz & 16 GB  \\ 
         \hline
     SQP-GS \cite{li2016sqp}    &  Dell
Precision T5810 &  4-core, 3.5 GHz  &  64 GB\\ 
         \hline
    \end{tabular}
    \egroup
    \vskip -1.5em
    \label{tab:machinea}
\end{table}

\section{Conclusions}\label{sec:conclusion}
This paper presents a novel convexified SSSC-OPF formulation based on a sufficient condition of small-signal stability, which is not based upon eigenvalue calculation and localized linearization of the dynamical model of power systems. The key feature of this convexification is the computational efficiency that allows for the scalability of SSSC-OPF. The solution feasibility is also recovered by incorporating a set of objective penalization functions. The BMI-based stability conditions are replaced by a more tractable trace-based upper bound in the form of a vector-norm. The proposed convexification technique has shown promising results through numerical simulations. The 9-bus, 39-bus and 118-bus systems' simulations show that the proposed method has achieved a stable optimal solution with a sufficiently low stability-induced cost. Future work will focus on improving the applicability by providing a guarantee of convergence and developing tighter convex relaxations.

\section*{Acknowledgement}
The authors are supported by NTU SUG, MOE AcRF TIER 1- 2019-T1-001-119 (RG 79/19), EMA \& NRF EMA-EP004-EKJGC-0003, and NRF DERMS for Energy Grid 2.0.

\bibliographystyle{IEEEtran}
\bibliography{main}

\appendices

\section{LMI Representation for Relaxation of $\mathbf{F}$}\label{app:1}
In this section, we introduce a lemma which facilitates the LMI representation of (\ref{eq:MBMI1}) and (\ref{eq:rlxM}). 
 
 \begin{lemma} \label{lemma}
If $X\in \mathbb{S}^+$, then the matrix inequality $Y^T Y - X \preceq 0$ is equivalent to:
    \begin{align}
      \begin{bmatrix}
      X & Y^T\\
      Y & \mathbf{I}
      \end{bmatrix} \succeq 0.
    \end{align}

Also, if $X\in \mathbb{S}^+$ and $U \succ 0$, then 
    \begin{align}
      \begin{bmatrix}
      X-Y^T Y & V^T\\
      V & U 
      \end{bmatrix} \succeq 0 \Longleftrightarrow 
  \begin{bmatrix}
      X & Y^T & V^T\\
      Y & \mathbf{I} & \mathbf{O}\\
      V & \mathbf{O} & U  
      \end{bmatrix} \succeq 0.
    \end{align}
\begin{proof}
The proof can be constructed using Schur's complement and Lemma 2.1 \cite{boyd2004convex} directly. We omit proof here.
\end{proof}
\end{lemma}
 
 By using the first part of Lemma \ref{lemma}, the LMI representation of (\ref{eq:MBMI1}) is:
 \begin{align}\label{eq:L1}
     \mathbf{L}_1&:= \begin{bmatrix}
     M & J^T+Z^T \\
     J+Z & \mathbf{I}
     \end{bmatrix}.
\end{align}

Similarly from the second part of Lemma \ref{lemma}, we obtain the LMI formulation of \eqref{eq:rlxM} as the following:
\begin{align}\label{eq:L22}
     \mathbf{L}_2 & := \begin{bmatrix}
     M &  Z^T & J^T\\
     Z & \mathbf{I} & \mathbf{O}\\
     J & \mathbf{O} & \mathbf{I}
     \end{bmatrix}.
\end{align}

With \eqref{eq:L1} and \eqref{eq:L22}, the constraints in (\ref{eq:MBMI1}) and (\ref{eq:rlxM}) can be cast as LMI and can be solved efficiently using SDP.


\section{Convex Envelopes for Park's Transform }\label{app:2}
The nonconvex relations (\ref{eq:parks}) are relaxed using the McCormick envelopes for bilinear terms \cite{mccormick1976computability}. By defining $u_i=\sin{\delta_i}$ and $v_i=\cos{\delta_i}$ the Park's transformation equations for $k \in \mathcal{G}(i)$ will be: 
\begin{equation}\label{eq:parks}
\begin{aligned}
    V_{di}=V_k\sin{(\delta_i-\theta_k)}=V_{x_k}\,u_i - V_{y_k} \, v_i\\
     V_{qi}=V_k\cos{(\delta_i-\theta_k)}=V_{x_k}\, v_i + V_{y_k} \, u_i
\end{aligned}
\end{equation}

The following convex envelopes for these bilinear terms have been used by a number of previous works such as \cite{hijazi2017convex} as follows: 
\begin{equation}\label{eq:mcc}
        \begin{aligned}
    \langle  ab\rangle^u_M \equiv \begin{cases} 
    (ab)^u &\leq a^u\,b + a\,b^l - a^u\,b^l \\
    (ab)^u &\leq  a^l\,b + a\,b^u - a^l\,b^u \\
    \end{cases}\\
    \langle  ab\rangle^l_M \equiv \begin{cases}
    (ab)^l &\geq a^u\,b + a\,b^u - a^u\,b^u \\
    (ab)^l &\geq a^l\,b + a\,b^l - a^l\,b^l
    \end{cases}
        \end{aligned}
\end{equation}

 We use $\langle \cdot \rangle^u_M $ and $\langle \cdot \rangle^l_M $ to denote upper and lower bounds for McCormick envelope of the bilinear terms as in (\ref{eq:mcc}).

Thus, the convex relaxations of (\ref{eq:parks}) become:
\begin{equation}\label{eq:rlx_parks}
\begin{aligned}
    V_{di} \leq\langle V_{x_i}\,u_i \rangle^u_M - \langle V_{y_i} \, v_i\rangle^l_M ,\\
    V_{di} \geq \langle V_{x_i}\,u_i \rangle^l_M - \langle V_{y_i} \, v_i\rangle^u_M ,\\
     V_{qi} \leq \langle V_{x_i}\, v_i\rangle^u_M + \langle V_{y_i} \, u_i\rangle^u_M ,\\
     V_{qi} \geq \langle V_{x_i}\, v_i\rangle^l_M + \langle V_{y_i} \, u_i\rangle^l_M .
\end{aligned}
\end{equation}

Two other quadratic relations need to be imposed to ensure that SDP solution is a feasible equilibrium point. The first one is the trigonometric equality between sine and cosine of the load angle. For this, we introduce the vector variables $\mathbf{U_u} \in \mathbb{R}^{n_g}$ and $\mathbf{U}_v \in \mathbb{R}^{n_g}$. Thus, the set of convex constraints for the trigonometric equality $u^2_i+v^2_i=1$ is: 
\begin{subequations}\label{eq:delta}
    \begin{align}
    \mathbf{U}_u+\mathbf{U}_v=\mathbf{1}\label{eq:delta1},\\
    {U}_{u_i}\geq u^2_i,\label{eq:delta2}\\
    \quad {U}_{v_i}\geq v^2_i . \label{eq:delta3}
    \end{align}
    \vspace{-5mm}
\end{subequations}

Here, (\ref{eq:delta2}) is the convex quadratic relaxation of the quadratic equalities between the lifting variable ($\mathbf{U}_u,\, \mathbf{U}_v$) and linear variables ($u_i,\, v_i$). Similarly, the convex relaxation of quadratic equality between the $d-q$ axis terminal voltages ($V_d,\,V_q$) and node voltages ($V_x,\, V_y$) can be given as: 
\begin{align}\label{eq:WnWdq}
            W_{dq_{i,i}}+W_{dq_{m,m}}=W_{_{k,k}}+W_{_{l,l}} ,
\end{align}
for $i=1 \dots n_g, \, m=i+n_g,\, k= \mathcal{G}(i)$, and $l=k+n_b$.

\subsubsection{Feasible solution recovery for Park's transformation}
For relaxations of trigonometric relations (\ref{eq:delta}), the objective penalty functions obtained for \eqref{eq:delta2} and \eqref{eq:delta3} $\forall~ i\in \mathcal{G}$ are
\begin{align}
    h_4(\mathbf{U}_u,u) & = \sum_i\big \{ {U}_{u_i}-2u_{i_o}u_i+u_{i_o}u_{i_o}\big \},\\
     h_5(\mathbf{U}_v,v) & =\sum_i\big \{ {U}_{v_i}-2v_{i_o}v_i+v_{i_o}v_{i_o}\big \}.
\end{align}

\end{document}